\documentclass{article}
\usepackage[english,russian,ukrainian]{babel}
\usepackage{latexsym,amssymb,amsmath}
\usepackage[cp1251]{inputenc}
\usepackage{amsfonts,amssymb}

\newcounter{theorem}
\newcounter{lemma}

\,
\,

\begin{document}
\large

\centerline{\textbf{Curvilinear integral theorem for $G$-monogenic mappings}}
\centerline{\textbf{in the algebra of complex quaternion}}
\vskip 2mm
\centerline{Tetyana KUZMENKO}
\vskip 2mm
\noindent\textbf{Abstract.} For $G$-monogenic mappings taking values in the algebra of complex quaternion we prove a curvilinear analogue of the Cauchy integral theorem in the case where a curve of integration lies on the boundary of a domain.
\vskip 2mm
\noindent{\textbf{AMS 2010:}  30G35.}
\large

\vskip 5mm
\textbf{Introduction}

Let $\mathbb{H(C)}$ be the quaternion algebra over the field of complex numbers
$\mathbb{C}$, whose basis consists of the unit $1$ of the algebra and of the elements
$I,J,K$
satisfying the multiplication rules:
$$I^2=J^2=K^2=-1,$$
$$\,IJ=-JI=K,\quad JK=-KJ=I,\quad KI=-IK=J.$$

In the algebra $\mathbb{H(C)}$ there exists another basis $\{e_1,e_2,e_3,e_4\}$ such that multiplication table in a new
basis can be represented as (see, e. g., \cite{Cartan})
$$
\begin{tabular}{c||c|c|c|c|}
$\cdot$ & $e_1$ & $e_2$ & $e_3$ & $e_4$\\
\hline
\hline
$e_1$ & $e_1$ & $0$ & $e_3$ & $0$\\
\hline
$e_2$ & $0$ & $e_2$ & $0$ & $e_4$\\
\hline
$e_3$ & $0$ & $e_3$ & $0$ & $e_1$\\
\hline
$e_4$ & $e_4$ & $0$ & $e_2$ & $0$\\
\hline
\end{tabular}\,\,.
$$
The unit of the algebra can be decomposed as $1=e_1+e_2$.

Let us consider the vectors
$$i_1=e_1+e_2, \quad i_2=a_1e_1+a_2e_2, \quad i_3=b_1e_1+b_2e_2,$$
$a_k,b_k\in\mathbb{C},\,k=1,2$,
which are linearly independent over the field of real numbers $\mathbb{R}$.
It means that the equality $\alpha_1i_1+\alpha_2i_2+\alpha_3i_3=0$ for $\alpha_1,\alpha_2,\alpha_3\in\mathbb{R}$ holds if and only if
 $\alpha_1=\alpha_2=\alpha_3=0$.

In the algebra $\mathbb{H(C)}$ we consider the linear span $E_3:=\{\zeta=xi_1+yi_2+zi_3:x,y,z\in\mathbb{R}\}$ generated by the
vectors $i_1,i_2,i_3$ over the field $\mathbb{R}$. A set $S\subset\mathbb{R}^3$ is associated with the set $S_\zeta:=
\{\zeta=xi_1+yi_2+zi_3:(x,y,z)\in S\}$ in $E_3$. We also note that a
 topological property of a set $S_\zeta$ in $E_3$ understand as the same
  topological property of the set $S$ in $\mathbb{R}^3$.
For example, we will say that a curve $\gamma_\zeta\subset E_3$ is homotopic to
 a point if $\gamma\subset\mathbb{R}^3$ is homotopic to a point, etc.


We say (see \cite{Shpakivskiy-Kuzmenko}) that a continuous mapping $\Phi:\Omega_\zeta\rightarrow\mathbb{H(C)}$ \big(or $\widehat{\Phi}:\Omega_\zeta\rightarrow\mathbb{H(C)}$\big) is \emph{right-$G$-monogenic}
\big(or \emph{left-$G$-monogenic}\big) in a domain
$\Omega_\zeta\subset E_3$\,,
if $\Phi$ \big(or $\widehat{\Phi}$\big) is differentiable in the sense
 of the G\^{a}teaux at every point of\, $\Omega_\zeta$\,, i.~e. for every $\zeta\in
 \Omega_\zeta$ there exists an element $\Phi'(\zeta)\in\mathbb{H(C)}$ \big(or $\widehat{\Phi}'(\zeta)\in\mathbb{H(C)}$\big) such that
$$\lim\limits_{\varepsilon\rightarrow 0+0}\Big(\Phi(\zeta+\varepsilon h)-\Phi(\zeta)\Big)\varepsilon^{-1}= h\Phi'(\zeta)\quad\forall\,h\in E_3$$
$$\Biggr(\text{or }\,\, \lim\limits_{\varepsilon\rightarrow 0+0}
\left(\widehat{\Phi}(\zeta+\varepsilon h)-\widehat{\Phi}(\zeta)\right)
\varepsilon^{-1}= \widehat{\Phi}'(\zeta)h\quad\forall\,h\in E_3\Biggr).$$

The Cauchy integral theorems for holomorphic functions of the complex variable are
   fundamental results of the classical complex analysis.
Analogues of these results are also important tools in the quaternionic analysis.

In the paper \cite{Shp-Kuzm-intteor} were established some analogues of classical
integral theorems of the theory of analytic functions of the complex variable: the surface and curvilinear Cauchy
integral theorems and the Cauchy integral formula. The Morera theorem was proved in the paper \cite{Shp-Kuz-umb}. Taylor's and Laurent's expansions of $G$-monogenic mappings are obtained in \cite{Kuzm-series}.

Namely, in the paper \cite{Shp-Kuzm-intteor} was proved a curvilinear analogue of the Cauchy integral theorem in the case where a curve of integration lies in a domain of $G$-monogeneity.

In the present paper we prove the curvilinear Cauchy integral theorem for $G$-monogenic mappings in the case where a curve of integration lies on the boundary of a domain of $G$-monogeneity.
	
\vskip 5mm
\textbf{The main result}

Let $\gamma$ be a Jordan rectifiable curve in $\mathbb{R}^3$. For a continuous
mapping $\Psi:\gamma_\zeta\rightarrow\mathbb{H(C)}$
 of the form
\begin{equation}\label{Phi-form}
\Psi(\zeta)=\sum\limits_{k=1}^{4}\Big(U_k(x,y,z)+iV_k(x,y,z)\Big) e_k,
\end{equation}
where $(x,y,z)\in\gamma$ and $U_k:\gamma\rightarrow\mathbb{R}$,
$V_k:\gamma\rightarrow\mathbb{R}$,
we define integrals along a Jordan rectifiable curve $\gamma_\zeta$ by
 the equalities
$$\int\limits_{\gamma_\zeta}{d\zeta\,\Psi(\zeta)}:=\sum\limits_{k=1}^{4}{e_k\int
\limits_{\gamma}U_k(x,y,z)dx}+\sum\limits_{k=1}^{4}{i_2e_k\int\limits_{\gamma}
U_k(x,y,z)dy}+$$
$$+\sum\limits_{k=1}^{4}{i_3e_k\int\limits_{\gamma}U_k(x,y,z)dz}+i\sum
\limits_{k=1}^{4}{e_k\int\limits_{\gamma}V_k(x,y,z)dx}+$$
$$+i\sum\limits_{k=1}^{4}{i_2e_k\int\limits_{\gamma}V_k(x,y,z)dy}+i\sum
\limits_{k=1}^{4}{i_3e_k\int\limits_{\gamma}V_k(x,y,z)dz}$$

\noindent and

$$\int\limits_{\gamma_\zeta}{\Psi(\zeta)\,d\zeta}:=\sum\limits_{k=1}^{4}
{e_k\int\limits_{\gamma}U_k(x,y,z)dx}+\sum\limits_{k=1}^{4}{e_ki_2\int
\limits_{\gamma}U_k(x,y,z)dy}+$$
$$+\sum\limits_{k=1}^{4}{e_ki_3\int\limits_{\gamma}U_k(x,y,z)dz}+i\sum
\limits_{k=1}^{4}{e_k\int\limits_{\gamma}V_k(x,y,z)dx}+$$
$$+i\sum\limits_{k=1}^{4}{e_ki_2\int\limits_{\gamma}V_k(x,y,z)dy}+i\sum
\limits_{k=1}^{4}{e_ki_3\int\limits_{\gamma}V_k(x,y,z)dz},$$
where $d\zeta:=dxi_1+dyi_2+dzi_3$.

In the paper \cite{Shp-Kuzm-intteor} for right-$G$-monogenic mappings was obtained the following analogue of the Cauchy integral theorem.

\vskip 2mm
\textbf{Theorem 1 \cite{Shp-Kuzm-intteor}.} \emph{Let $\Phi:\Omega_\zeta\rightarrow\mathbb{H(C)}$ be a
right-$G$-monogenic mapping in a domain $\Omega_\zeta$. Then for every
closed Jordan rectifiable curve $\gamma_\zeta$ homotopic to a point in
$\Omega_\zeta$\,, the following equality is true:
\begin{equation}\label{int=0}
\int\limits_{\gamma_\zeta}d\zeta\,\Phi(\zeta)=0.
\end{equation}}
\vskip 2mm

Below we establish sufficient conditions for the curve $\gamma_\zeta$\, lying on the boundary $\partial\Omega_\zeta$ of a domain $\Omega_\zeta$\, such that the equality (\ref{int=0}) holds. For this goal we apply a scheme of the paper \cite{Pl-Shp-dop} for $G$-monogenic mappings.

Let on a boundary $\partial\Omega_\zeta$ of the domain $\Omega_\zeta$ given closed Jordan rectifiable curve $\gamma_\zeta\equiv\gamma_\zeta(t),$ where $0\leq t\leq1,$ homotopic to an interior point $\zeta_0\in\Omega_\zeta$. It means that there exists the mapping $H(s,t)$ continuous on the square $[0,1]\times[0,1]$, such that $H(0,t)=\gamma_\zeta(t),\, H(1,t)\equiv\zeta_0$, and all curves $\gamma_\zeta^s\equiv\gamma_\zeta^s(t):=\{\zeta=H(s,t):0\leq t\leq1\}$ for $0<s<1$ are contained in the domain $\Omega_\zeta$.

Consider also the curves $\Gamma_\zeta^t\equiv\Gamma_\zeta^t(s):=\{\zeta=H(s,t):0\leq s\leq1\}$. Denote by $\Gamma[\zeta_1,\zeta_2]$ arc of Jordan oriented rectifiable curve, beginning at the point $\zeta_1$ and ending at the point $\zeta_2$, and denote by the mes a linear Lebesgue measure of a rectifiable curve.

As in the paper \cite{Shp-Kuz-umb}, for the element $\zeta=xi_1+yi_2+zi_3$ we define
the Euclidian norm
$$\|\zeta\|=\sqrt{x^2+y^2+z^2}.$$
Using the Theorem of equivalents of norms, for the element
$a:=\sum\limits_{k=1}^{4}(a_{1k}+ia_{2k})e_k$,\, $a_{1k},a_{2k}\in\mathbb{R}$, we have
the following inequalities
$$|a_{1k}+ia_{2k}|\leq\sqrt{\sum\limits_{k=1}^{4}\big(a_{1k}^2+a_{2k}^2\big)}\,\leq c \|a\|,$$
where $c$ is a positive constant does not dependent on $a$.

\vskip 2mm
\textbf{Theorem 2.} \emph{Suppose that $\Phi:\overline{\Omega}_\zeta\rightarrow\mathbb{H(C)}$ is a continuous  mapping in the closure $\overline{\Omega}_\zeta$ of a domain $\Omega_\zeta$ and right-$G$-monogenic in $\Omega_\zeta$. Suppose also that $\gamma_\zeta\subset\partial\Omega_\zeta$ is a closed Jordan rectifiable curve homotopic to an interior point $\zeta_0\in\Omega_\zeta$\,, the curves of the family $\{\Gamma_\zeta^t:0\leq t\leq1\}$ are rectifiable and the set} $\{\text{mes}\,\gamma_\zeta^s:0\leq\ s\leq1\}$ \emph{is bounded, then the equality \em(\ref{int=0})\em\, is true.}
\vskip 2mm

\textbf{Proof.} Let $\varepsilon>0$. We fix the number $\rho\in\big(0,\frac{1}{2}\,\text{mes}\,\gamma_\zeta\big)$ such that for arbitrary $\zeta_1,\zeta_2\in\overline{\Omega}_\zeta$ from the condition $||\zeta_1-\zeta_2||<2\rho$ follows the inequality
\begin{equation}\label{teor-8-1}
||\Phi(\zeta_1)-\Phi(\zeta_2)||<\varepsilon.
\end{equation}

Since the mapping $H$ is uniformly continuous on the square $[0,1]\times[0,1]$, then there exists $\delta>0$ such that for all $s\in(0,\delta)$ and $t,t'\in[0,1]:|t-t'|<\delta$ the inequality $|H(0,t)-H(s,t')|<\rho$ is true.

Let numbers $0=t_0<t_1<\ldots<t_n<1$ such that for corresponding points $\zeta_{0,k}:=H(0,t_k)$ of the curve $\gamma_\zeta$\, the following relations are fulfilled
$$\text{mes}\,\gamma_\zeta[\zeta_{0,k},\zeta_{0,k+1}]=\rho \quad \text{for} \quad k=\overline{0,n-1},$$
$$\text{mes}\,\gamma_\zeta[\zeta_{0,n},\zeta_{0,0}]\leq\rho.$$

It is obvious that $2\leq n\leq\left[\frac{\text{mes}\,\gamma_\zeta}{\rho}\right]+1$.

Let us consider the points $\zeta_{s,k}:=H(s,t_k)$ of the curve $\gamma_\zeta^s$ and the curves
$$\Upsilon_{[k]}^s:=\gamma_\zeta[\zeta_{0,k},\zeta_{0,k+1}]\cup\Gamma_\zeta^{t_{k+1}}[\zeta_{0,k+1},\zeta_{s,k+1}]\cup\gamma_\zeta^s[\zeta_{s,k+1},\zeta_{s,k}]\cup\Gamma_\zeta^{t_k}[\zeta_{s,k},\zeta_{0,k}]$$
for $k=\overline{0,n},$ where $\zeta_{s,n+1}:=\zeta_{s,0}$ for $0\leq s\leq1$, setting that the orientation of curves $\Upsilon_{[k]}^s$ is induced by orientation of the curve $\gamma_\zeta.$

Let $s\in(0,\delta).$ Since for all $\zeta\in\Upsilon_{[k]}^s$ the inequality $||\zeta-\zeta_{0,k}||\leq2\rho$ is true, then by Theorem 2 \cite{Shp-Kuzm-intteor}, Lemma 4.1 \cite{Shp-Kuz-umb} and the inequality (\ref{teor-8-1}), we have
$$\Bigg\|\int\limits_{\gamma_\zeta}d\zeta\Phi(\zeta)\Bigg\|=\Bigg\|\sum_{k=0}^n\int\limits_{\Upsilon_{[k]}^s}d\zeta(\Phi(\zeta)-\Phi(\zeta_{0,k}))\Bigg\|\leq$$
$$\leq c\sum\limits_{k=0}^n\int\limits_{\Upsilon_{[k]}^s}||d\zeta||\,||\Phi(\zeta)-\Phi(\zeta_{0,k})||\leq c\varepsilon\sum\limits_{k=0}^n\text{mes}\,\Upsilon_{[k]}^s\leq$$
$$\leq c\varepsilon\Big(\text{mes}\,\gamma_\zeta+\text{mes}\,\gamma_\zeta^s+2(n+1)\max\limits_{k=\overline{0,n}}\,\text{mes}\,\Gamma_\zeta^{t_k}[\zeta_{s,k},\zeta_{0,k}]\Big)\leq$$
\begin{equation}\label{teor-8-2}
\leq M\varepsilon\Bigg(1+\frac{1}{\rho}\max\limits_{k=\overline{0,n}}\,\text{mes}\,\Gamma_\zeta^{t_k}[\zeta_{s,k},\zeta_{0,k}]\Bigg),
\end{equation}
and a constant $M$ does not depend on $\varepsilon$ and $\rho$.

Passing to the limit in the inequality (\ref{teor-8-2}) as $s\rightarrow0$, we have the inequality
$$\Bigg\|\int\limits_{\gamma_\zeta}d\zeta\Phi(\zeta)\Bigg\|\leq M\varepsilon.$$
Now passing to the limit in the last inequality as $\varepsilon\rightarrow0$, we obtain the equality (\ref{int=0}). The Theorem is proved.
\vskip 2mm

The similar statement is true for the left-$G$-monogenic mappings.

\vskip 2mm
\textbf{Theorem 3.} \emph{Suppose that $\widehat{\Phi}:\overline{\Omega}_\zeta\rightarrow\mathbb{H(C)}$ is a continuous mapping in the closure $\overline{\Omega}_\zeta$ of a domain $\Omega_\zeta$ and left-$G$-monogenic in $\Omega_\zeta$. Suppose also that $\gamma_\zeta\subset\partial\Omega_\zeta$ is a closed Jordan rectifiable curve homotopic to an interior point $\zeta_0\in\Omega_\zeta$\,, the curves of the family $\{\Gamma_\zeta^t:0\leq t\leq1\}$ are rectifiable and the set} $\{\text{mes}\,\gamma_\zeta^s:0\leq\ s\leq1\}$ \emph{is bounded, then the following equality is true:}
$$\int\limits_{\gamma_\zeta}\widehat{\Phi}(\zeta)d\zeta=0.$$
\vskip 2mm

\renewcommand{\refname}{References}

\end{document}